\documentclass[11pt]{amsart}

\usepackage{graphicx}

\newtheorem{proposition}{Proposition}
\newtheorem{corollary}[proposition]{Corollary}
\newtheorem{question}{Question}

\def\tb{tb}
\def\mindeg{\operatorname{min-deg}}
\def\maxdeg{\operatorname{max-deg}}
\def\breadth{\operatorname{breadth}}
\def\tb{\operatorname{tb}}
\def\sl{\operatorname{sl}}
\def\maxtb{\overline{\tb}}
\def\maxsl{\overline{\sl}}

\begin{document}

\title{On Arc Index and Maximal Thurston--Bennequin Number}
\author{Lenhard Ng}
\address{Mathematics Department, Duke University, Durham, NC 27708}
\email{ng@math.duke.edu}
\urladdr{http://alum.mit.edu/www/ng}

\begin{abstract}
We discuss the relation between arc index, maximal
Thurston--Bennequin number, and Khovanov homology for knots.
As a consequence, we calculate the arc index and maximal
Thurston--Bennequin number for all knots with at most $11$ crossings.
For some of these knots, the calculation requires a consideration of
cables which also allows us to compute the maximal self-linking
number for all knots with at most $11$ crossings.
\end{abstract}

\maketitle

\section{Introduction and Results}
\label{sec:intro}

Let $K$ be a knot in $S^3$. Define a \textit{grid diagram} of $K$ to
be an oriented knot
diagram for $K$ consisting of a union of horizontal and vertical line
segments, such that at every crossing, the vertical segment crosses
over the horizontal segment. Any knot has a grid diagram.
In the literature, grid diagrams or their
equivalents have gone by many alternate names, including ``arc
presentations'', ``asterisk presentations'',
``square-bridge presentations'', and ``fences''. Grid diagrams have
been much studied lately, most recently because of their use in the
combinatorial definition of knot Floer homology \cite{bib:MOS}; for
background on grid diagrams, see, e.g., \cite{bib:Cro}.

The \textit{arc number} of a grid diagram is the number of horizontal
(or, equivalently, vertical) segments in the diagram.
The \textit{arc index} of $K$, written $\alpha(K)$,
is the minimal arc number over all grid diagrams for $K$.

It is well-known that grid diagrams are closely related to Legendrian
knots from contact geometry (see, e.g., \cite{bib:Et} for an introduction to
Legendrian knots and a more geometric description of the invariants
$\tb$ and $\sl$ below). A front for a Legendrian knot can be obtained by 
rotating any grid diagram slightly counterclockwise and eliminating
each corner by either smoothing it out or replacing it by a
cusp. Conversely, any Legendrian knot can be represented by a grid
diagram.

In this context, the \textit{Thurston--Bennequin number} $\tb$ and
\textit{self-linking number} $\sl$ of a grid
diagram $G$, which are invariants of the associated Legendrian knot,
can be defined as follows. Let $w(G)$ denote the writhe of 
$G$; let $c(G)$
denote the number of lower-right, ``southeast'',
corners of $G$ (these correspond to
the right cusps of the Legendrian front); and let $c_{\downarrow}(G)$
denote the number of southeast corners oriented down and to the left,
plus the number of northwest corners oriented to the left and down
(these correspond to the downward-oriented cusps of the Legendrian front).
Then
\begin{align*}
\tb(G) &= w(G) -
c(G) \\
\sl(G) &= w(G) - c_{\downarrow}(G).
\end{align*}
We remark that the self-linking number is usually defined for
transverse rather than Legendrian knots; $\sl$ defined here is the
self-linking number of the positive transverse pushoff of the
Legendrian knot, and can be expressed as $\tb(G)-\operatorname{r}(G)$,
where $\operatorname{r}(G)$ is the rotation number of the Legendrian
knot. 

The \textit{maximal Thurston--Bennequin number} of a knot $K$,
written $\overline{\tb}(K)$,
is the maximal $\tb$ over all grid diagrams for $K$; similarly, the
\textit{maximal self-linking number} $\overline{\sl}(K)$ is the
maximal $\sl$ over all grid diagrams for $K$. It is not hard to see that
$\overline{\tb}(K) \leq \overline{\sl}(K)$ for all $K$, while it is an important classical result of Bennequin \cite{bib:Ben} that $\overline{\sl}(K) < \infty$ for any $K$.
Calculating $\maxtb$ and $\maxsl$ is of natural interest to knot
theorists, particularly since each provides a lower bound for various
topological knot invariants, including the slice genus $g_4$
\cite{bib:Rud1} and the concordance invariants $\tau$ \cite{bib:Plam1}
and $s$ \cite{bib:Plam2,bib:Shu}.

There is a fundamental relation between arc index and the maximal
Thurston--Bennequin numbers of a knot $K$ and its mirror
$\overline{K}$, first described by Matsuda in \cite{bib:Mat}:
\begin{equation}
- \alpha(K) \leq \maxtb(K) + \maxtb(\overline{K}).
\label{eq:arc}
\end{equation}
The proof of this
inequality is short and we recall it here. Consider a grid diagram for
$K$ with arc number $\alpha(K)$. This diagram produces a Legendrian knot
of topological type $K$, as described above, as well as a Legendrian
knot of type $\overline{K}$, by rotating the diagram slightly less
than $90^{\circ}$ clockwise, changing every crossing, and smoothing
the corners. Then it is easy to see that the Thurston--Bennequin
numbers of these two Legendrian knots sum to $-\alpha(K)$.

Equation (\ref{eq:arc}) leads to an approach to calculate arc index
and maximal Thurston--Bennequin number for specific knots, as
follows:

\begin{enumerate}
\renewcommand{\labelenumi}{(\alph{enumi})}
\item
find a possibly minimal grid diagram of $K$;
\item
find upper bounds for $\maxtb(K)$ and $\maxtb(\overline{K})$
individually, or for their sum;
\item
see if equality is forced to hold in (\ref{eq:arc}).
\end{enumerate}

\noindent This approach (essentially) has been used to calculate arc
index for alternating knots \cite{bib:BP} and knots with up to $10$
crossings \cite{bib:Bel}. In both cases, the upper bound in step (b)
is provided by the Kauffman polynomial.

In this note, we apply this approach to knots with at most $11$
crossings, using grid diagrams provided by Baldwin and Gillam
\cite{bib:BG} and the Khovanov bound for $\maxtb$ \cite{bib:NgKho}.
We compute arc index and maximal Thurston--Bennequin number
for all knots with at most $11$ crossings.
Let $\mindeg$ and $\maxdeg$ denote the minimum
and maximum degrees of a Laurent polynomial in the specified
variable, let $\breadth = \maxdeg - \mindeg$, and let $Kh_K(q,t)$
denote the two-variable Poincar\'e polynomial for $\mathfrak{sl}_2$
Khovanov homology.

\begin{proposition}
Let $K$ be a knot with $11$ or fewer crossings. We have
\[
\alpha(K) = \breadth_q Kh_K(q,t/q)
\]
with the following exceptions: $\alpha(10_{124}) = 8$,
$\alpha(10_{132}) = 9$,
$\alpha(11n_{12}) = 10$, $\alpha(11n_{19}) = 9$,
$\alpha(11n_{38}) = 9$, $\alpha(11n_{57}) = 10$,
$\alpha(11n_{88}) = 10$, and $\alpha(11n_{92}) = 10$. Here the
chirality of the knot is irrelevant.
\label{prop:arc}
\end{proposition}

\begin{proposition}
Let $K$ be a knot with $11$ or fewer crossings. We have
\[
\maxtb(K) = \mindeg_q Kh_K(q,t/q)
\]
with the following exceptions:
\label{prop:tb}
\begin{align*}
\maxtb (\overline{10_{124}}) &= -15 & \maxtb(\overline{11n_{38}}) &=
-4 \\
\maxtb (\overline{10_{132}}) &= -1 & \maxtb (\overline{11n_{57}}) &=-13\\
\maxtb (11n_{12}) &= -2 & \maxtb (\overline{11n_{88}}) &=-13 \\
\maxtb (11n_{19}) &= -8 & \maxtb (11n_{92}) &=-6.
\end{align*}
\end{proposition}

\noindent
The $\maxtb$ data from Proposition~\ref{prop:tb} for knots with up to
$11$ crossings can be found online at \texttt{KnotInfo} \cite{bib:KnotInfo}.

The exceptional cases in Proposition~\ref{prop:tb} require
strengthening previously known upper bounds for $\maxtb$ and are
presented in Section~\ref{sec:proofs}.
The computation of $\maxtb$ for $11n_{19}$
uses a strengthening of the Kauffman bound on $\maxtb$ derived
from work of Rutherford \cite{bib:Ru} and a subsequent observation of
K\'alm\'an \cite{bib:Kal}; the computation of $\maxtb$ for
$\overline{10_{132}}$, $11n_{12}$, $\overline{11n_{38}}$,
$\overline{11n_{57}}$, $\overline{11n_{88}}$, and $11n_{92}$ uses
cable links.

Nutt \cite{bib:Nutt} previously directly computed arc index for all
knots with $9$ or fewer crossings, and Beltrami \cite{bib:Bel}, as
mentioned earlier, extended this computation to knots with $10$
crossings. The author \cite{bib:NgKho} previously computed maximal
Thurston--Bennequin number for all knots with $10$ or fewer
crossings except $\overline{10_{132}}$.

Josh Greene \cite{bib:Gr} has proposed the following very
interesting question:

\begin{question}
Does a grid diagram realizing the arc index of a knot necessarily
realize the maximal Thurston--Bennequin number for the knot? An
equivalent statement is that
\begin{equation}
- \alpha(K) = \maxtb(K) + \maxtb(\overline{K})
\label{eq:q1}
\end{equation}
for all knots $K$.
\label{q1}
\end{question}

No counterexamples are currently known. In particular, we have the
following consequence of Propositions~\ref{prop:arc} and
\ref{prop:tb}:

\begin{corollary}
(\ref{eq:q1}) holds for all knots $K$ with $11$ or fewer crossings.
\end{corollary}

Greene notes that (\ref{eq:q1}) also holds for alternating knots
by \cite{bib:BP} and the fact that the Kauffman bound for $\maxtb$ is
sharp for alternating knots \cite{bib:NgKho,bib:Ru}, and for torus
knots by Etnyre and Honda's classification of Legendrian torus knots
\cite{bib:EH}.

We conclude this section with a discussion of maximal self-linking number.
There is an intriguing analogy between $\maxtb$ and $\maxsl$:
\[
\text{arc index } : \text{ braid index } :: ~\maxtb~:~\maxsl.
\]
Keiko Kawamuro \cite[Conjecture 3.2]{bib:Kaw} has made a conjecture
which can be restated as follows to
parallel Question~\ref{q1}:

\begin{question}
Does a braid whose closure is a particular knot,
with a minimal number of strands (the
braid index), necessarily realize the maximal self-linking number for
the knot? An equivalent statement is that
\begin{equation}
-2b(K) = \maxsl(K) + \maxsl(\overline{K})
\label{eq:q2}
\end{equation}
for all knots $K$, where $b(K)$ is the braid index of $K$.
\label{q2}
\end{question}

\noindent Note that (\ref{eq:q2}), like (\ref{eq:q1}), holds if $=$
is replaced by $\leq$. The celebrated MFW inequality
\cite{bib:FW,bib:Mor} gives a lower bound for braid index and an upper
bound for $\maxsl$ in terms of the HOMFLY-PT polynomial $P_K(a,z)$:
\[
-2b(K) \leq \maxsl(K) + \maxsl(\overline{K}) \leq
-\breadth_a P_K(a,z) - 2.
\]
Thus the answer to Question~\ref{q2} is ``yes'' for all knots for
which the ``weak'' MFW inequality $2b(K) \geq \breadth_a P_K(a,z) + 2$
is sharp.

In fact, more is true.
In Section~\ref{ssec:proofs2}, we calculate $\maxsl(K)$ for the
$5$ knots with at most $10$ crossings for which MFW is not sharp.
This calculation, combined with an analogous calculation by T.\
Khandhawit \cite{bib:Khan} for the $14$ knots with $11$ crossings where MFW
is not sharp, yields the following result.

\begin{proposition}
Let $K$ be a knot with $11$ or fewer crossings. We have
\[
\maxsl(K) = -\maxdeg_a P_K(a,z) - 1
\]
with the following exceptions:
\label{prop:sl}
\begin{align*}
\maxsl(\overline{9_{42}}) &= -5 &
\maxsl(\overline{11n_{24}}) &= -5 &
\maxsl(11n_{86}) &= -3 \\
\maxsl(\overline{9_{49}}) &= -11 &
\maxsl(\overline{11n_{33}}) &= -7 &
\maxsl(\overline{11n_{117}}) &= -7 \\
\maxsl(\overline{10_{132}}) &= -1 &
\maxsl(\overline{11n_{37}}) &= -3 &
\maxsl(\overline{11n_{124}}) &= -7 \\
\maxsl(\overline{10_{150}}) &= -9 &
\maxsl(\overline{11n_{70}}) &= -7 &
\maxsl(\overline{11n_{136}}) &= -13 \\
\maxsl(10_{156}) &= -7 &
\maxsl(\overline{11n_{79}}) &= -7 &
\maxsl(\overline{11n_{171}}) &= -13 \\
\maxsl(11n_{20}) &= -7 &
\maxsl(\overline{11n_{82}}) &= -5 &
\maxsl(\overline{11n_{180}}) &= -13 \\
&&&&
\maxsl(\overline{11n_{181}}) &= -13.
\end{align*}
\end{proposition}

\begin{corollary}
(\ref{eq:q2}) holds for all knots $K$ with $11$ or fewer crossings.
\end{corollary}

As for (\ref{eq:q1}), no counterexamples to (\ref{eq:q2}) are currently
known.

\section*{Acknowledgments}

I am grateful to Dror Bar-Natan, Danny Gillam, Jeremy Green, Josh
Greene, Matt Hedden, Tamas K\'alm\'an, Keiko Kawamuro, and Tirasan
Khandhawit for useful
conversations, and to Princeton University, the
University of Virginia, and Tom Mark for their hospitality during the
course of this work. This work was partially supported by NSF grant
DMS-0706777.

\section{Proofs}
\label{sec:proofs}

In this section, we provide more details for the discussion in
Section~\ref{sec:intro}, and prove the main
results. Section~\ref{ssec:proofs1} proves Proposition~\ref{prop:arc},
and Proposition~\ref{prop:tb} for all but six cases; Section~\ref{ssec:proofs2}
uses cables to fill in the remaining cases and also prove
Proposition~\ref{prop:sl}.

\subsection{Arc index and $\maxtb$}
\label{ssec:proofs1}

Two very useful bounds for $\maxtb$ are the \textit{Kauffman bound}
\cite{bib:FT,bib:Rud1,bib:Tab}
\begin{equation}
\maxtb(K) \leq - \maxdeg_a F_K(a,z)-1,
\label{eq:Kauffman}
\end{equation}
where $F_K$ is the two-variable Kauffman polynomial of $K$, and
the \textit{Khovanov bound} \cite{bib:NgKho}
\begin{equation}
\maxtb(K) \leq \mindeg_q Kh_K(q,t/q),
\label{eq:Kh}
\end{equation}
where $Kh_K$ is the Poincar\'e polynomial for $\mathfrak{sl}_2$
Khovanov homology.\footnote{Note: There are many different conventions
  regarding knot chirality in the literature. These results, and this
  paper in general, use the
  conventions that conform to the Knot Atlas \cite{bib:BN}. In
  particular, the Kauffman bound (\ref{eq:Kauffman})
  uses the opposite convention for the Kauffman polynomial to the
  one used in many Legendrian-knot papers, including
  \cite{bib:FT,bib:Ng2bridge,bib:Tab}.}
It was noted in \cite{bib:NgKho} that the Khovanov bound is at least
as strong as the Kauffman bound for all knots with $11$ or fewer
crossings, although the two bounds are incommensurate in general.

Combining (\ref{eq:arc}) and (\ref{eq:Kauffman}) yields
\begin{equation}
\alpha(K) \geq \breadth_a F_K(a,z) +2.
\label{eq:arcKauffman}
\end{equation}
The inequality (\ref{eq:arcKauffman}) is originally due to Morton
and Beltrami \cite{bib:MB}, and Beltrami \cite{bib:Bel} used it to
compute the arc index of all $10$-crossing knots. Bae and Park
\cite{bib:BP} proved that (\ref{eq:arcKauffman}) is sharp (i.e.,
equality holds) for alternating knots, where both sides are equal to
the crossing number plus $2$.

Combining (\ref{eq:arc}) and (\ref{eq:Kh}) instead yields the
following result.

\begin{proposition}
If $K$ is a knot, then
\begin{equation}
\alpha(K) \geq \breadth_q Kh_K(q,t/q).
\label{eq:arcKh}
\end{equation}
If $K$ has a grid diagram with arc number equal to $\breadth_q
Kh_K(q,t/q)$, then (\ref{eq:arcKh}) is sharp, as is the Khovanov bound
for both $\maxtb(K)$ and $\maxtb(\overline{K})$.
\label{prop:arcKh}
\end{proposition}

We now apply Proposition~\ref{prop:arcKh} to prove
Proposition~\ref{prop:arc}.

\begin{proof}[Proof of Proposition~\ref{prop:arc}]
Because of the behavior of arc index and Khovanov homology under
connected sum, it suffices to consider prime knots only. In
addition, the result holds for alternating knots $K$; here
$\alpha(K) = \breadth_a F_K(a,z) +2 = c(K) + 2$, where $c(K)$ is the
crossing number of $K$, and both Kauffman and Khovanov bounds for
$\maxtb$ are sharp \cite{bib:NgKho,bib:Ru}.

Baldwin and Gillam \cite{bib:BG}, with the help of the program
\texttt{Gridlink} \cite{bib:Grid}, have constructed grid diagrams
for all nonalternating prime knots with $11$ or fewer crossings;
these presentations, which include a few diagrams constructed by the
author, are available at
\verb+http://www.math.brown.edu/~wgillam/hfk/+. For most of these
diagrams, the arc number is equal to $\breadth_q Kh_K(q,t/q)$, as
can easily be checked by computer. (The author used
\texttt{KnotTheory} \cite{bib:BN} for this computation.) The exceptions are
$10_{124}$, $10_{132}$, $11n_{12}$, $11n_{19}$, $11n_{38}$, $11n_{57}$,
$11n_{88}$, and $11n_{92}$; for each of these, however, arc index has
been computed in \cite{bib:Nutt}.
\end{proof}

Before proving Proposition~\ref{prop:tb}, we introduce a minor
strengthening of the Kauffman bound (\ref{eq:Kauffman}),
due to K\'alm\'an \cite{bib:Kal} and based on work of Rutherford
\cite{bib:Ru}. Rutherford's paper relates the Dubrovnik version of
the Kauffman polynomial, $D_K(a,z) = F_K(ia,-iz)$, to certain
partitions of fronts of Legendrian knots known as rulings
\cite{bib:ChP}.

\begin{proposition}[K\'alm\'an]
Let $K$ be a knot, and let $p_K(z)$ denote the polynomial in $z$
which is the leading term of $F_K(ia,-iz)$ with respect to $a$.
\label{prop:improved}
If $p_K(z)$ does not have all nonnegative coefficients, then
\[
\maxtb(K) \leq - \maxdeg_a F_K(a,z)-2.
\]
\end{proposition}

\begin{proof}
Suppose that the Kauffman bound (\ref{eq:Kauffman}) is sharp for
$K$, and consider a Legendrian knot $L$ of type $K$ for which
$\tb(L) = -\maxdeg_a F_K(a,z)-1$. By \cite{bib:Ru}, we have $p_K(z)
=\sum_{\rho\in\Gamma(L)} z^{j(\rho)}$, where $\Gamma(L)$ is the set
of rulings of $L$ and $j$ is an integer-valued function on rulings.
In particular, $p_K(z)$ has all nonnegative coefficients.
\end{proof}

Proposition~\ref{prop:improved} allows us to lower the Kauffman
bound by $1$ in some cases. Unfortunately, it does not apply to many
small knots. The hypotheses of the proposition apply to seven knots
with $11$ crossings or fewer: $\overline{10_{136}}$, $11n_{19}$,
$11n_{20}$, $\overline{11n_{37}}$, $\overline{11n_{50}}$,
$11n_{86}$, and $\overline{11n_{126}}$. For six of these, the
improved Kauffman bound is only as good as the Khovanov bound
(\ref{eq:Kh}); for $11n_{19}$, however, it improves on the Khovanov
bound as well, to yield $\maxtb(11n_{19}) \leq -8$. For
$12$-crossing knots, Proposition~\ref{prop:improved} yields the best
known bound on $\maxtb$ for three knots, according to the tabulation
from \texttt{KnotInfo} \cite{bib:KnotInfo}: $\maxtb(12n_{25})\leq
-5$, $\maxtb(\overline{12n_{502}}) \leq -17$,
$\maxtb(\overline{12n_{603}}) \leq -12$.

One can similarly use Rutherford's work to obtain an improved
HOMFLY-PT bound on $\maxtb$, when the leading coefficient of the
HOMFLY-PT polynomial does not have all nonnegative coefficients, and a
``mixed'' improved bound when the HOMFLY-PT and Kauffman bounds agree
and the leading coefficient of their difference does not have
all nonnegative coefficients. These seem to be applicable to fewer cases
than the improved Kauffman bound, however.

We can now prove Proposition~\ref{prop:tb}.

%
%
\begin{proof}[Proof of Proposition~\ref{prop:tb}]
As in the proof of Proposition~\ref{prop:arc}, the result holds
unless $K$ is one of the knots $10_{124}$, $10_{132}$, $11n_{12}$, $11n_{19}$,
$11n_{38}$, $11n_{57}$, $11n_{88}$, or $11n_{92}$, with either
chirality.

As discussed earlier, the case $10_{124} = T(3,5)$ is covered by
\cite{bib:EH}; $\maxtb(10_{124}) = 7$ and $\maxtb(\overline{10_{124}})
= -15$, and the Khovanov bound is sharp for the former but not for the latter.
For $11n_{19}$, Proposition~\ref{prop:improved} gives
$\maxtb(11n_{19}) \leq -8$, while both Kauffman and Khovanov bounds
give $\maxtb(\overline{11n_{19}}) \leq -1$; since
$\alpha(11n_{19})=9$ by Nutt's table \cite{bib:Nutt}, these bounds
for $\maxtb(11n_{19})$ and $\maxtb(\overline{11n_{19}})$ are
sharp.

The remaining cases, $10_{132}$, $11n_{12}$, $11n_{38}$, $11n_{57}$,
$11n_{88}$, and $11n_{92}$, are addressed by Corollary~\ref{cor:cable}
in the next section. (In fact,
$10_{124}$ and $11n_{19}$ can also be addressed in the same way.)
\end{proof}

\subsection{Cables, $\maxtb$, and $\maxsl$}

\label{ssec:proofs2}

Suppose that we wish to assemble a table of maximal
Thurston--Bennequin and self-linking numbers for small knots.
There are several knots with $11$ or fewer crossings for which
all of the known general upper bounds on $\maxtb$ or $\maxsl$ fail to
be sharp: $7$ for $\maxtb$, $19$ for $\maxsl$. What can one do in these cases?
One case for $\maxtb$, $\overline{10_{124}}$, is the $(3,-5)$ torus knot, and the
classification of Legendrian torus knots due to Etnyre and Honda
\cite{bib:EH} shows that $\maxtb(\overline{10_{124}}) = -15$; the best
general upper bound gives $\maxtb(\overline{10_{124}}) \leq -14$. For
the other cases, however, there is no classification result. For
these, we turn to cable links.

If $K$ is a knot, let $D_n(K)$ denote the $n$-framed double
($2$-cable link) of $K$, where both components of $D_n(K)$ are
oriented the same way as $K$. Our strategy is to bound $\maxtb$ and
$\maxsl$ for $D_n(K)$ from above via one of the standard bounds, and
then use these upper bounds to bound $\maxtb$ and $\maxsl$ for $K$ via
the following easy result.\footnote{The observation that (\ref{eq:sldouble})
holds for all $n$, not only $n=0$, is due to Khandhawit
\cite{bib:Khan}.}

\begin{proposition}
We have
\begin{equation}
\maxtb(D_n(K)) \geq \begin{cases} 2 \,\maxtb(K) + 2n, & n> \maxtb(K)
  \\
4n, & n\leq \maxtb(K)
\end{cases}
\label{eq:tbdouble}
\end{equation}
and
\begin{equation}
\maxsl(D_n(K)) \geq 2\,\maxsl(K) + 2n.
\label{eq:sldouble}
\end{equation}
As a consequence of (\ref{eq:tbdouble}), if $\maxtb(D_n(K)) < 2m+2n$
for some $m,n$ with $m\leq n$, then $\maxtb(K) < m$.
\label{prop:cable}
\end{proposition}

\begin{proof}
We first prove (\ref{eq:tbdouble}).
Let $L$ be a Legendrian knot of type $K$. Define the ``Legendrian double''
$D(L)$ to be the Legendrian link whose front is given by two copies of
$L$ offset slightly in the vertical ($z$) direction; then $D(L)$ is
topologically the $\tb(L)$-framed double of $K$, and $\tb(D(L)) =
4\tb(L)$.

If $n \leq \maxtb(K)$, then choose $L$ such that $\tb(L) = n$. Since
$D(L)$ is topologically $D_n(K)$ and $\tb(D(L)) = 4n$, it follows that
$\maxtb(D_n(K)) \geq 4n$. If $n > \maxtb(K)$, then choose $L$ such
that $\tb(L) = \maxtb(K)$. Add $n-\maxtb(K)$ positive twists to the
framing on $D(L)$ by inserting $n-\maxtb(K)$ pieces of the form
\includegraphics[height=12pt]{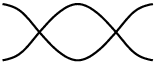} into the front
of $D(L)$ to obtain a Legendrian link $D'(L)$ which is topologically
$D_n(K)$. Each of the pieces adds $2$ to $\tb$, and so $\tb(D'(L)) =
4\,\maxtb(K)+2(n-\maxtb(K))$; it follows that $\maxtb(D_n(K)) \geq
2\,\maxtb(K)+2n$.

To prove (\ref{eq:sldouble}),
we use the alternate formulation, first observed by Bennequin
\cite{bib:Ben}, for self-linking number in terms 
of braids. If $B$ is a braid of $m$ strands and writhe (algebraic
crossing number) $w$, then define $\sl(B) = w-m$; $\maxsl(K)$ is
the maximum value of $\sl(B)$ over all braids $B$ whose closure is
$K$.

Given $K$, let $B$ be a braid whose closure is $K$ for which
$\sl(B) = \maxsl(K)$. Construct a double $B'$ of $B$ with $2m$ strands
consisting of two slightly offset copies of $B$; in algebraic terms,
replace each generator $\sigma_i^{\pm 1}$ in the braid word for $B$ by
$(\sigma_{2i}\sigma_{2i-1}\sigma_{2i+1}\sigma_{2i})^{\pm 1}$ to obtain
$B'$. If $w$ is the writhe of $B$, then the closure of $B'$ is
$D_w(K)$, and $\sl(B') = 4w-2m$.
Add in $n-w$ positive twists to the beginning of $B'$ (i.e.,
append $\sigma_1^{2n-2w}$ to the braid word for $B'$) to obtain another
braid $B''$ with $2m$ strands. The closure of $B''$ is $D_n(K)$, and
$\sl(B'') = 2w-2m+2n$.
It follows that $\maxsl(D_n(K)) \geq 2w-2m+2n = 2\,\maxsl(K)+2n$.
\end{proof}

\begin{corollary}
The values of $\maxtb$ for $\overline{10_{132}}$,
$11n_{12}$, $\overline{11n_{38}}$,
$\overline{11n_{57}}$, $\overline{11n_{88}}$, and $11n_{92}$ (and
their mirrors) are as given in Proposition~\ref{prop:tb}.
\label{cor:cable}
\end{corollary}

\begin{proof}
We combine the Khovanov bound for $\maxtb(D_n(K))$ with
Proposition~\ref{prop:cable}. For instance, the Khovanov bound yields
$\maxtb(D_3(\overline{10_{132}})) \leq 5$, which with
Proposition~\ref{prop:cable} implies that $\maxtb(\overline{10_{132}})
\leq -1$. The Khovanov bound also shows directly that
$\maxtb(10_{132}) \leq -8$; from Proposition~\ref{prop:arcKh} and
(\ref{eq:arc}), we conclude that $\maxtb(10_{132}) = -8$ and
$\maxtb(\overline{10_{132}}) = -1$.

Similarly, the Khovanov bound gives $\maxtb(D_3(11n_{12})) \leq 3$,
$\maxtb(D_1(\overline{11n_{38}})) \leq -6$,
$\maxtb(D_{-7}(\overline{11n_{57}})) \leq -39$,
$\maxtb(D_{-7}(\overline{11n_{88}})) \leq -39$, and
$\maxtb(D_{-1}(11n_{92}) \leq -13$, and these bounds produce the
values of $\maxtb$ for $11n_{12}$, $\overline{11n_{38}}$,
$\overline{11n_{57}}$, $\overline{11n_{88}}$, and $11n_{92}$ given in
Proposition~\ref{prop:tb}. We remark that these doubles are links with
$40+$ crossings, and computing their Khovanov homology is not
altogether trivial. The particular framings of the doubles were
chosen to try to minimize crossings, and each Khovanov homology was
computed using the program \texttt{JavaKh}, written by Jeremy
Green, within \texttt{KnotTheory} \cite{bib:BN}.
\end{proof}

%
%
%

We next use Proposition~\ref{prop:cable} to prove
Proposition~\ref{prop:sl}. 

\begin{proof}[Proof of Proposition~\ref{prop:sl}]
We prove the result for knots with $10$ or fewer crossings, and refer
the reader to \cite{bib:Khan} for the $11$-crossing case, which is
proved by the same technique.
There is nothing to prove if the weak MFW inequality $2b(K) \geq
\breadth_a P_K(a,z) + 2$ is sharp. There are five knots with $10$ or
fewer crossings for which equality does not hold for MFW: $9_{42}$,
$9_{49}$, $10_{132}$, $10_{150}$, and $10_{156}$. For these, we use
the HOMFLY-PT bound on $\maxsl(D_0(K))$ and
Proposition~\ref{prop:cable} to bound $\maxsl$. The HOMFLY-PT bound
yields an upper bound on $\maxsl(D_0(K))$ of $-8$, $-20$, $0$, $-16$, and
$-12$, respectively. (For some of these computations, the author found
the program \texttt{K2K} \cite{bib:K2K} to be useful.)
These give the exceptional values for $\maxsl$ in
the statement of Proposition~\ref{prop:sl}.

For example, since $\maxsl(D_0(\overline{9_{42}})) \leq -8$,
Proposition~\ref{prop:cable} implies that $\maxsl(\overline{9_{42}})
\leq -4$; since the self-linking number for any knot is odd, it
follows that $\maxsl(\overline{9_{42}}) \leq -5$. The usual HOMFLY-PT
bound also implies that $\maxsl(9_{42}) \leq -3$. Since $b(9_{42}) =
4$ and $-2b(9_{42}) \leq \maxsl(9_{42}) + \maxsl(\overline{9_{42}})$,
equality holds everywhere.
\end{proof}

We close with two remarks. First, using cables along the lines
presented here is not entirely new; Stoimenow \cite{bib:Sto} showed
that $\overline{10_{132}}$ is not quasipositive using almost identical
methods.

Second, in the situations where the general upper
bounds for $\maxtb(K)$ and $\maxsl(K)$ (Kauffman, Khovanov, HOMFLY-PT) fail
to be sharp, it seems that one can often apply these bounds to the
double or perhaps general $m$-cable of $K$ to deduce a sharp bound for
$\maxtb(K)$ and $\maxsl(K)$. Proposition~\ref{prop:cable} has
a straightforward analogue for $m$-component cables of $K$.
 For instance, if $C_m(K)$ denotes the
$0$-framed $m$-component cable of $K$, then
\[
\maxsl(C_m(K)) \geq m\,\maxsl(K).
\]
It seems at least within the realm of possibility that
$\maxsl(K) = \lim_{m\to\infty} \maxsl(C_m(K))/m$, and that the
HOMFLY-PT bound for $\maxsl(C_m(K))$ might in general give a sharp
bound for $\maxsl(K)$ for \textit{all} $K$. A similar but slightly more
complicated statement could hold for $\maxtb$.

Thus there might be a way to calculate $\maxtb$ and $\maxsl$ for all
knots, by applying the general upper bounds to cables. We note,
however, that calculating these upper bounds for cables is generally
quite computationally intensive and may be infeasible for
``medium-sized'' knots of, say, $12$ crossings or more.


\end{document}